\documentclass{amsart}
\usepackage{amssymb,latexsym}
\usepackage{mltex}
\usepackage{amsmath}
\usepackage{epsfig}
\usepackage{graphics}
 \usepackage{lscape} 
\usepackage{enumerate}
\usepackage{amsthm}
\usepackage{hyperref}

\newcommand{\be}{\begin{enumerate}}  \newcommand{\ee}{\end{enumerate}}
\newcommand{\beqt}{\begin{equation}}  \newcommand{\eeqt}{\end{equation}}
\newcommand{\beq}{\begin{eqnarray}}  \newcommand{\eeq}{\end{eqnarray}}
\newcommand{\beQ}{\begin{eqnarray*}} \newcommand{\eeQ}{\end{eqnarray*}}

\newcommand{\RR}{{\mathbb R}}

\newcommand{\SSS}{{\mathbb S}}

\newcommand{\tr}{\mathrm{tr}}

\newcommand{\id}{\mathrm{Id}}

\newcommand{\lgra}{\longrightarrow}

\def\Id{\operatorname{Id}}

\newtheorem{example}{Exemples}[section]
\newtheorem{thm}{Theorem}[section]
\newtheorem{lemma}[thm]{Lemma}
\newtheorem{prop}[thm]{Proposition}
\newtheorem{cor}[thm]{Corollary}
\newtheorem{remark}[thm]{Remark}
\newtheorem{remarks}[thm]{Remarques}
\newtheorem{definition}[thm]{Definition}
\newtheorem{notation}[thm]{Notation}
\newtheorem{exabout:ample}[thm]{Exemple}

\begin{document}

\title{A Fundamental Theorem for submanifolds of multiproducts of real space forms}
\author{Marie-Am\'elie Lawn and Julien Roth}
%\address{Julien Roth, Laboratoire d'Analyse et de Math\'ematiques Appliqu\'ees, UPEM-UPEC, CNRS, F-77454 Marne-la-Vall\'ee}
%\email{julien.roth@u-pem.fr}
\subjclass[2010]{53A42, 53C20, 53C21}
\keywords{Almost-umbilical hypersurfaces, space forms, constant scalar curvature}

\maketitle

\begin{abstract}
We prove a Bonnet theorem for isometric immersions of submanifolds into the products of an arbitrary number of simply connected real space forms. Then, we prove the existence of associated families of minimal surfaces in such products. Finally, in the case of $\SSS^2\times\SSS^2$, we give a {\it complex} version of the main theorem in terms of the two canonical complex structures of  $\SSS^2\times\SSS^2$.
\end{abstract}

%\tableofcontents

\section{Introduction}
It is a classical problem of submanifold theory  to determine when a
Riemannian manifold $(M^n,g)$ can be immersed into a fixed Riemannian
manifold $(\bar{M}^{n+p}, \bar{g})$.  The well-known Gauss, Ricci and
Codazzi equations relate the intrinsic and extrinsic curvatures, and
any submanifold of any Riemannian manifold must satisfy them.
Conversely, the classical Bonnet theorem \cite{Bo} states that on a
surface, given first and second fundamental forms satisfying the Gauss
and Codazzi equations, this surface is locally embeddable into the
Euclidean $3$-space $\RR^3$. This result can be generalized to higher
codimension \cite{Ten}, and the classical  Fundamental Theorem of
Submanifolds states that, in fact, the Gauss, Codazzi and Ricci
equations are necessary and sufficient conditions for a Riemannian
$n$-dimensional manifold to admit a (local) immersion into a space of
constant sectional curvature of dimension $n+d$. \\
If the ambient space is not of constant sectional curvature, proving
fundamental theorems is technically difficult and there are few
results known. Moreover, the Gauss, Codazzi and Ricci equations are in
general not sufficient anymore and other conditions are required in
order to produce the immersion. In \cite{Da2}, Daniel gave such a
characterization for surfaces in the three-dimensional Thurston
geometries with four-dimensional isometry groups, by computing the
Christoffel symbols explicitly and using the technique of Cartan
moving frames. In higher dimensions, he also stated in \cite{Da}
necessary and sufficient conditions for an $n$-dimensional Riemannian
manifold to be isometrically immersible into the products
$\mathbb{S}^n\times\mathbb{R}$ and $\mathbb{H}^n\times\mathbb{R}$,
also using the moving frame technique. This allowed him to study the
existence of associated families in the case of minimal surfaces. This
result was later generalized by the second author \cite{Ro2} in the
case where the ambient space is a Lorentzian product. Very recently,
Ortega and the first author \cite{LO} proved fundamental theorems
characterizing immersions of hypersurfaces into (quasi-)Einstein
manifolds, specifically Robertson-Walker warped products. These spaces
play an important role in standard models of cosmology, arising as
solutions of the  non-vacuum Einstein equations, and have therefore a
great importance in Lorentzian geometry. As an application, conditions
were obtained for $3$-dimensional hypersurfaces in Robertson-Walker
spacetimes to be foliated by surfaces whose mean curvature vector is
either lightlike or zero (including maximal surfaces, marginally outer
trapped surfaces (MOTS), and mixed cases), hence providing an helpful
tool for the study of horizons on Robertson-Walker spacetimes with
spacelike or timelike causal character, including marginally outer
trapped tubes. \\
Extending the result of \cite{Da}, Kowalczyk and
Lira-Tojeiro-Vit\'orio proved independently in \cite{Ko} and
\cite{LTV} the existence and uniqueness of isometric immersions in a
product of two spaces forms of constant sectional curvature. In this
paper we generalize their result  to immersions into multiproducts
$\tilde{P}=M_{1}\times\dots\times M_m$ of real space form of arbitrary
dimension and arbitrary sectional curvature. The key idea is to use
the projections $\pi_i$, $i\in{1,\dots ,m}$ into each of the factors
of the product. Each projection induces then two operators on the
tangent bundle and two operators on the normal bundle of the
submanifold satisfying some properties and some compatibility
equations which can be deduced from the Gauss and Weingarten formulas
of the immersion. We prove that, conversely, these conditions together
with the Gauss, Codazzi and Ricci equations are necessary and
sufficient conditions to immerse a Riemannian manifold isometrically
into such an ambient space.\\
As an application we then prove the existence of a one-parameter
associated family of isometric immersions for minimal surfaces in
multiproducts. Finally we consider the special case where the ambient
space is $\mathbb{S}^2\times\mathbb{S}^2$ and give a complex version
of our fundamental theorem in terms of the induced complex structures.\\
The authors want to thank F. Torralbo for helpful discussions.

\section{Multiproducts of space forms and their submanifolds}\label{sec2}
We consider the product space $(\widetilde{P}=M_1\times \cdots\times M_m, \widetilde{g}=g_1\oplus\cdots\oplus g_m)$ where $(M_i,g_i)$ is the simply connected real space form of dimension $n_i$ and constant sectional curvature $c_i$. Moreover, without loss of generality, we assume that $c_i\neq0$ for $i\in\{1,\cdots,m-1\}$ and that $c_m$ may possibly be zero. We denote by $\pi_i$ the projection of any tangent vector $X$ on $TM_k$. These projections satisfy the following relations
$$\displaystyle\left\{
\begin{array}{l}
\widetilde{g}(\pi_iX,Y)=\widetilde{g}(X,\pi_iY)\ \text{for any}\ X,Y\in\Gamma(TP),\\
\pi_i\circ\pi_i=\pi_i,\\
\pi_i\circ\pi_j=0\ \text{if}\ i\neq j,\\
\widetilde{\nabla}\pi_i=0,\\
rank(\pi_i)=n_i,\\

\sum_{i=1}^m\pi_i={\rm Id}_{T\widetilde{P}}.
\end{array}
\right.$$
Moreover, the curvature tensor of $\widetilde{P}$ is given by
\begin{eqnarray}\label{extcurvature}\widetilde{R}(X,Y)Z=\sum_{i=1}^mc_i\bigl[\langle\pi_i Y,\pi_iZ\rangle\pi_iX-\langle\pi_i X,\pi_iZ\rangle\pi_iY\bigr].\end{eqnarray}
Now, we consider a Riemannian manifold $(M^n,g)$ isometrically immersed into $\widetilde{P}$. We denote by $NM$ the normal bundle and by $\nabla^{\perp}$ the normal connection and by $B:TM\times TM\lgra NM$ the second fundamental form. For any $\xi\in NM$, $A_{\xi}$ is the Weingarten operator associated to $\xi$ and defined by $\widetilde{g}(A_{\xi}X,Y)=\widetilde{g}(B(X,Y),\xi)$, with $X,Y$ vectors tangent to $M$.\\
 For any $i\in\{1,\cdots,m\}$, the projection $\pi_i$ induces the existence of the following four operators
$f_i:TM\rightarrow TM$, $h_i:TM\rightarrow NM$, $s_i:NM\rightarrow TM$ and $t_i:NM\rightarrow NM$, such that
\begin{eqnarray}\label{relationfhst}
\pi_iX=f_iX+h_iX\quad\text{and}\quad
\pi_i\xi=s_i\xi+t_i\xi.
\end{eqnarray} 
From the symmetry of  the $\pi_i$, we obtain that, for any $i\in\{1,\cdots,m\}$, $f_i$ and $t_i$ are symmetric and for any $X\in\Gamma(TM)$ and $\xi\in\Gamma(NM)$
\beqt\label{relation1.0}
\widetilde{g}(h_iX,\xi)=\widetilde{g}(X,s_i\xi).
\eeqt
 In addition, from the fact that $\displaystyle\sum_{i=1}^m\pi_i=Id_{T\widetilde{P}}$, we get the following identities
\beqt\label{relation1.0} 
\sum_{i=1}^mf_i=Id_{TM},\quad \sum_{i=1}^mt_i=Id_{E}, \quad \sum_{i=1}^ms_i=0\quad\text{and}\quad \sum_{i=1}^mh_i=0.
\eeqt
Moreover, we have the following relations between these operators coming from the fact that $\pi_i\circ\pi_j=\delta_i^j\pi_i$
\begin{align}
&f_i\circ f_j+s_i\circ h_j=\delta_i^jf_i,& \label{relation1.1}\\
&t_i\circ t_j+h_i\circ s_j=\delta_i^jt_i,& \label{relation1.2}\\
&f_i\circ s_j+s_i\circ t_j=\delta_i^js_i,&\label{relation1.3}\\
&h_i\circ f_j+t_i\circ h_j=\delta_i^jh_i,&\label{relation1.4}
\end{align}
where $\delta_i^j$ is the classical Kronecker symbol, that is, $1$ if $i=j$ and $0$ if $i\neq j$. 
Moreover, from the fact that $\pi_i$ is parallel, we deduce easily that for any $X,Y\in TM$ and $\xi\in NM$, we have
\begin{align}
&\nabla_X(f_iY)-f_i(\nabla_XY)=A_{h_iY}X+s_i(B(X,Y)),& \label{relation2.1}\\
&\nabla^{\perp}_X(h_iY)-h_i(\nabla_XY)=t_i(B(X,Y))-B(X,f_iY),& \label{relation2.2}\\
&\nabla_X^{\perp}(t_i\xi)-t_i(\nabla^{\perp}_X\xi)=-B(s_i\xi,X)-h_i(A_{\xi}X),& \label{relation2.3}\\
&\nabla_X(s_i\xi)-s_i(\nabla^{\perp}_X\xi)=-f_i(A_{\xi}X)+A_{t_i\xi}X.&\label{relation2.4}
\end{align}
Finally, from the expression for the curvature tensor $\widetilde{R}$, we get the following Gauss, Codazzi and Ricci equations
\begin{equation}\label{Gauss}\tag{G}
R(X,Y)Z=\sum_{i=1}^mc_i\bigg[ \left\langle f_iY,Z\right\rangle f_iX-\left\langle f_iX,Z\right\rangle f_iY \bigg]+A_{B(Y,Z)}X-A_{B(X,Z)}Y,
\end{equation}
\begin{equation}\label{Codazzi}\tag{C}
(\nabla_XB)(Y,Z)-(\nabla_YB)(X,Z)=\sum_{i=1}^mc_i\bigg[ \left\langle f_iY,Z\right\rangle h_iX-\left\langle f_iX,Z \right\rangle h_iY\bigg],
\end{equation}
\begin{equation}\label{Ricci}\tag{R}
R^{\perp}(X,Y)\xi=\sum_{i=1}^mc_i\bigg[ \left\langle h_iY,\xi\right\rangle h_iX-\left\langle h_iX,\xi\right\rangle h_iY\bigg]+B(A_{\xi}Y,X)-B(A_{\xi}X,Y).
\end{equation}
\section{Main result}
Now, conversely, consider $(M^n,g)$ a Riemannian manifold and $E$ a $d$-dimensional vector bundle over $M$ endowed with a metric $\overline{g}$ and a compatible connection $\overline{\nabla}$. Moreover, let  $B:TM\times TM\lgra E$ be a symmetric $(2,1)$-tensor and $f_i:TM\lgra TM$, $h_i:TM\lgra E$ and $t_i:E\lgra E$ be some $(1,1)$-tensors for $i\in\{1,\cdots,m\}$. We define $s_i$ as the dual of $h_i$ with respect to the metrics $\widetilde{g}:=g\oplus\overline{g}$ on $TM\oplus E$, that is, for any $X\in T_xM$ and $\xi\in E_x$,
$$\overline{g}_x(h_iX,\xi)=g_x(X,s_i\xi).$$
Finally, for any $\xi\in\Gamma(E)$, we define $A_{\xi}$ by 
$$\langle A_{\xi}X,Y\rangle=\langle B(X,Y),\xi\rangle,$$
for any $X,Y\in\Gamma(TM)$.
Following the discussions of Section \ref{sec2}, we introduce now the following natural definition.
\begin{definition}\label{defcomp}
We say that $(M,g,E,\overline{g},\overline{\nabla},B,f_i,h_i,t_i)$ satisfies the compatibility equations for the multiproduct $\widetilde{P}=M_1\times \cdots\times M_m$ if 
\begin{enumerate}[i)]
\item $f_i$ and $t_i$ are symmetric for any $i\in\{1,\cdots,m\}$,
\item for any $i\in\{1,\cdots,m\}$, Equations \eqref{relation1.0}-\eqref{relation2.3} are satisfied, that is,
\begin{align}
&\sum_{i=1}^mf_i=Id_{TM},\quad \sum_{i=1}^mt_i=Id_{E}, \quad \sum_{i=1}^ms_i=0\quad\text{and}\quad \sum_{i=1}^mh_i=0&\nonumber\\
&f_i\circ f_j+s_i\circ h_j=\delta_i^jf_i,&\nonumber\\
&t_i\circ t_j+h_i\circ s_j=\delta_i^jt_i,&\nonumber \\
&f_i\circ s_j+s_i\circ t_j=\delta_i^js_i,&\nonumber\\
&h_i\circ f_j+t_i\circ h_j=\delta_i^jh_i,&\nonumber\\
&\nabla_X(f_iY)-f_i(\nabla_XY)=A_{h_iY}X+s_i(B(X,Y)),& \nonumber\\
&\overline{\nabla}_X(h_iY)-h_i(\nabla_XY)=t_i(B(X,Y))-B(X,f_iY),& \nonumber\\
&\overline{\nabla}_X(t_i\xi)-t_i(\overline{\nabla}_X\xi)=-B(s_i\xi,X)-h_i(A_{\xi}X),&\nonumber
%&\nabla_X(s_i\xi)-s_i(\nabla^{\perp}_X\xi)=-f_i(A_{\xi}X)+A_{t_i\xi}X.&\nonumber
\end{align}
\item the rank of $\pi_i$ is $n_i$ for any $i\in\{1,\cdots,m\}$ and $\sum_{i=1}^mn_i=n+p,$, and
\item the Gauss, Ricci and Codazzi equations \eqref{Gauss}, \eqref{Codazzi} and \eqref{Ricci} are satisfied. Namely for any $X,Y,Z\in\Gamma(TM)$ and any $\xi\in\Gamma(E)$,
\begin{align}
&R(X,Y)Z=\sum_{i=1}^mc_i\bigg[ \left\langle f_iY,Z\right\rangle f_iX-\left\langle f_iX,Z\right\rangle f_iY \bigg]+A_{B(Y,Z)}X-A_{B(X,Z)}Y,&\nonumber\\
&(\nabla_XB)(Y,Z)-(\nabla_YB)(X,Z)=\sum_{i=1}^mc_i\bigg[ \left\langle f_iY,Z\right\rangle h_iX-\left\langle f_iX,Z \right\rangle h_iY\bigg],&\nonumber\\
&\overline{R}(X,Y)\xi=\sum_{i=1}^mc_i\bigg[ \left\langle h_iY,\xi\right\rangle h_iX-\left\langle h_iX,\xi\right\rangle h_iY\bigg]+B(A_{\xi}Y,X)-B(A_{\xi}X,Y),\nonumber&
\end{align}
where $\overline{R}$ is the curvature associated with the connection $\overline{\nabla}$.
\end{enumerate}
\end{definition}
We can now state the main result of the paper.
\begin{thm}\label{thm1}
Let $(M^n,g)$ be a simply connected  Riemannian manifold and $E$ a $d$-dimensional vector bundle over $M$ endowed with a metric $\overline{g}$ and a compatible connection $\nabla^{\perp}$. Moreover, let  $B:TM\times TM\lgra E$ be a symmetric $(2,1)$-tensor and $f_i:TM\lgra TM$, $h_i:TM\lgra E$ and $t_i:E\lgra E$ be some $(1,1)$-tensors for $i\in\{1,\cdots,m\}$. If $(M,g,E,\overline{g},\overline{\nabla},B,f_i,h_i,t_i)$ satisfies the compatibility equations for the multiproduct $\widetilde{P}=M_1\times \cdots\times M_m$ then, there exists an isometric immersion $\varphi:M\lgra\widetilde{P}$ such that the normal bundle of $M$ for this immersion is isomorphic to $E$ and such that the second fundamental form, the normal connection and the projections on each factor of $T\widetilde{P}_{|M}$ are given by $B$ , $\overline{\nabla}$ and $(f_i,h_i,t_i)$ respectively. Precisely, there exists a vector bundle isometry $\widetilde{\varphi}:E\lgra (\varphi(M))^{\perp}$ so that
$$\pi_i(\varphi_*X)=\varphi_*(f_iX)+\Phi(h_iX),$$
$$\pi_i(\widetilde{\varphi}\xi)=\varphi_*(s_iX)+\Phi(t_i\xi),$$
$$II_f=\widetilde{\varphi}\circ B,$$
$$\nabla^{\perp}\widetilde{\varphi}=\widetilde{\varphi}\overline{\nabla}.$$
Moreover, this isometric immersion is unique up to an isometry of $\widetilde{P}$.
\end{thm}
Our approach to prove this theorem is not based on the moving frame technique, but is in the spirit of \cite{Ko} and uses techniques introduced in \cite{Di} and \cite{DNV}.\\ \\
{\bf Proof:} 
We give the proof for the case $c_m\neq0$, the case $c_m=0$ can be proved analogously with minor changes. First, for any $i\in\{1,\dots,m\}$, let us denote by $E_i$ a trivial line bundle over $M$ equipped with the Euclidean metric, if $c_i>0$, and minus the Euclidean metric, if $c_i<0$. We consider the vector bundle $F$ over M,
$$ F=TM\oplus E \overset{m}{\underset{i=1}{\oplus}}E_i,$$
defined by the orthogonal Withney sum of Riemannian vector bundles. We denote by $\widetilde{g}$ the metric over $F$ obtained from $g$, $\overline{g}$ and the metrics on each $E_i$. For any $i\in\{1,\dots,m\}$, we consider a section $\xi_i$ of $E_i$ such that $\widetilde{g}(\xi_i,\xi_i)=\frac{1}{c_i}$. We introduce now the following connection on $F$, denoted by $D$
\begin{align}
&D_XY=\nabla_XY+B(X,Y)-\sum_{i=1}^mc_ig(f_iX,Y)\xi_i,&\nonumber\\
&D_X\nu=\nabla^{\perp}_X\nu-A_{\nu}X-\sum_{i=1}^mc_i\overline{g}(h_iX,\nu)\xi_i,&\nonumber\\
&D_X\xi_i=f_iX+h_iX,\nonumber
\end{align}
for any vector fields $X,Y$ tangent to $M$ and any section $\nu$ of $E$.
\begin{lemma}
The connection $D$ is compatible with the metric $\widetilde{g}$.
\end{lemma}
{\bf Proof:} This comes easily from the definition. Let $X,Y,Z\in\Gamma(TM)$, $\nu,\eta\in\Gamma(E)$. We have
\beQ
X\widetilde{g}(Y,Z)&=&Xg(Y,Z)\\
&=&{g}(\nabla_XY,Z)+g(Y,\nabla_XZ)\\
&=&\widetilde{g}(D_XY,Z)+\widetilde{g}(Y,D_XZ),
\eeQ
since the tangential parts of $D_XY$ and $D_XZ$ are $\nabla_XY$ and $\nabla_XZ$ respectively. Similarly, we have
\beQ
X\widetilde{g}(\nu,\eta)&=&X\overline{g}(\nu,\eta)\\
&=&\overline{g}(\nabla^{\perp}_X\nu,\eta)+\overline{g}(\nu,\nabla^{\perp}_X\eta)\\
&=&\widetilde{g}(D_X\nu,\eta)+\widetilde{g}(\nu,D_X\eta),
\eeQ
since the normal parts of $D_X\nu$ and $D_X\eta$ are $\nabla^{\perp}_X\nu$ and $\nabla^{\perp}_X\eta$ respectively. Moreover, we have
\beQ
X\widetilde{g}(\xi_i,\xi_j)&=0\\
&=&\widetilde{g}(f_iX+h_iX,\xi_j)+\widetilde{g}(\xi_i,f_jX+h_jX)\\
&=&\widetilde{g}(D_X\xi_i,\xi_j)+\widetilde{g}(\xi_i,D_X\xi_j).
\eeQ
Finally for mixed terms, we have
\beQ
X\widetilde{g}(\xi_i,Y)&=0\\
&=&g(f_iX,Y)-c_i\widetilde{g}(\xi_i,\xi_i)g(f_iX,Y)\\
&=&\widetilde{g}(f_iX+h_iX,Y)+\widetilde{g}\left(\xi_i,\nabla_XY+B(X,Y)-\sum_{i=1}^mc_ig(f_iX,Y)\xi_i\right)\\
&=&\widetilde{g}(D_X\xi_i,\xi_j)+\widetilde{g}(\xi_i,D_X\xi_j),
\eeQ
and
\beQ
X\widetilde{g}(\xi_i,\nu)&=0\\
&=&g(h_iX,\nu)-c_i\widetilde{g}(\xi_i,\xi_i)g(h_iX,\nu)\\
&=&\widetilde{g}(f_iX+h_iX,\nu)+\widetilde{g}\left(\xi_i,\nabla^{\perp}_X\nu-A_{\nu}X-\sum_{i=1}^mc_ig(h_iX,\nu)\xi_i\right)\\
&=&\widetilde{g}(D_X\xi_i,\nu)+\widetilde{g}(\xi_i,D_X\nu).
\eeQ
By bilinearity, we get the property for any sections $\alpha$ and $\beta$ of $F$.
\hfill$\square$ \\ \\
Now,  we consider the curvature tensor associated with the connection $D$, denoted by $\mathcal{R}^D$ and defined by $\mathcal{R}^D(X,Y)=D_XD_Y-D_YD_X-D_{[X,Y]}$. We can prove the following
\begin{lemma}
The connection $D$ is flat, that is, $\mathcal{R}^D=0$.
\end{lemma}
{\bf Proof:} Let $X,Y,Z\in\Gamma(TM)$ and $\nu\in\Gamma(E)$. We will prove that $\mathcal{R}^D(X,Y)Z=0$, $\mathcal{R}^D(X,Y)\nu=0$ and $\mathcal{R}^D(X,Y)\xi_i=0$ for any $i\in\{1,\cdots,m\}$. Then by linearity of the curvature  $\mathcal{R}^D$ in its third argument, we will get that $\mathcal{R}^D=0$. First, we have
\beQ
D_XD_YZ&=&\nabla_X\nabla_YZ+B(X,\nabla_YZ)-\sum_{i=1}^mc_ig(f_iX,\nabla_YZ)\xi_i+\nabla^{\perp}_XB(Y,Z)-A_{B(Y,Z)}X\\
&&-\sum_{i=1}^mc_i\Bigl[\overline{g}(B(Y,Z),h_iX)\xi_i+(f_iX+h_iX)+g(\nabla_Xf_iY,Z)\xi_i+g(f_iY,\nabla_XZ)\xi_i\Bigr].
\eeQ
Therefore, we get
\beQ
\mathcal{R}^D(X,Y)Z&=&R(X,Y)Z-\sum_{i=1}^mc_i\Bigl[g(f_iY,Z)f_iX-g(f_iX,Z)f_iY\Bigl]-A_{B(Y,Z)}X+A_{B(X,Z)}Y\\
&&+(\nabla_X B)(Y,Z)-(\nabla_YB)(X,Z)-\sum_{i=1}^mc_i\Bigl[g(f_iY,Z)h_iX-g(f_iX,Z)h_iY\Bigl]\\
&&-\sum_{i=1}^mc_i\Bigl[g\big((\nabla_Xf_i)Y-A_{h_iY}X+s_iB(X,Y),Z\big)\Bigr]\\
&&-\sum_{i=1}^mc_i\Bigl[g\big((\nabla_Yf_i)X-A_{h_iX}Y+s_iB(X,Y),Z\big)\Bigr]\\
&=&0
\eeQ
by using the Gauss equation (first line), the Codazzi equation (second line) and equation \eqref{relation2.1} (third and fourth lines). Similarly, we have
\beQ
D_XD_YZ\nu&=&\nabla^{\perp}_X\nabla^{\perp}_Y\nu-\nabla_XA_{\nu}Y-B(X,A_{\nu}Y)-A_{\nabla^{\perp}_X\nu}Y\\
&&+\sum_{i=1}^mc_i\Bigl[g(f_iX,A_{\nu}Y)\xi_i+\overline{g}(\nu,h_iY)(f_iX+h_iX)\Bigl]\\
&&-\sum_{i=1}^mc_i\Bigl[\overline{g}(\nabla_Y^{\perp}\nu,h_iX)-\overline{g}(\nabla_X^{\perp}\nu,h_iY)\overline{g}(\nu,\nabla_X^{\perp}h_iY)\Bigl].
\eeQ
And hence,
\beQ
\mathcal{R}^D(X,Y)\nu&=&R^{\perp}(X,Y)\nu-B(X,A_{\nu}Y)+B(Y,A_{\nu}X)-\sum_{i=1}^mc_i\Bigl[\overline{g}(h_iY,\nu)h_iX-\overline{g}(h_iX,\nu)h_iY\Bigl]\\
&&+\nabla_YA_{\nu}X+A_{\nabla_Y^{\perp}\nu}X-\nabla_XA_{\nu}Y-A_{\nabla_X^{\perp}\nu}Y-\sum_{i=1}^mc_i\Bigl[\overline{g}(h_iY,\nu)f_iX-\overline{g}(h_iX,\nu)f_iY\Bigl]\\
&&+\sum_{i=1}^mc_i\Bigl[g(f_iX,A_{\nu}Y)-\overline{g}(\nu,(\nabla_Yh_i)X)+\overline{g}(t_iB(X,Y),\nu)\Bigr]\\
&&+\sum_{i=1}^mc_i\Bigl[g(f_iY,A_{\nu}X)-\overline{g}(\nu,(\nabla_Xh_i)Y)+\overline{g}(t_iB(X,Y),\nu)\Bigr].
\eeQ
The first line in the right hand side vanishes due to Ricci equation. The second line vanishes by Codazzi equation and the third and fourth lines vanish by using equation \eqref{relation2.2}.
Finally, for any $i\in\{1,\cdots,m\}$, we have
\beQ
D_XD_Y\xi_j&=&\nabla_Xf_jY+B(X,f_jY)-\sum_{i=1}^mc_ig(f_iX,f_iY)\xi_i\\
&&+\nabla_X^{\perp}h_jY-A_{h_jY}X-\sum_{i=1}^mc_ig(h_iX,h_iY)\xi_i.\
\eeQ
Hence,
\beQ
\mathcal{R}^D(X,Y)\xi_j&=&(\nabla_Xf_j)Y-A_{h_jY}X-s_jB(X,Y)-(\nabla_Yf_j)X+A_{h_jX}Y+s_jB(X,Y)\\
&&\nabla_X^{\perp}h_jY+B(X,f_jY)-t_jB(X,Y)-\nabla_Y^{\perp}h_jX-B(Y,f_jX)-t_jB(X,Y)\\
&=&0
\eeQ
by equations \eqref{relation2.1} and \eqref{relation2.2}. Thus, we get that the connection $D$ is flat.
\hfill$\square$\\ \\
We define now for any $i\in\{1,\cdots,m\}$ the map $\pi_i:F\lgra F$ by  
$$\pi_iX=f_iX+h_iX,$$
$$\pi_i\nu=s_i\nu+t_i\nu,$$
$$\pi_i\xi_j=\delta_i^j\xi_j,$$
for any $X\in\Gamma(TM)$ and any $\nu\in\Gamma(E)$. We have the following properties
\begin{lemma}\label{lempi}
For any $i\in\{1,\cdots,m\}$, the map $\pi_i$ is symmetric with respect to $\widetilde{g}$, parallel with respect to $D$ and 
\begin{enumerate}
\item $\pi_i\circ\pi_i=\pi_i$,
\item $\sum_{i=1}^m\pi_i=Id_F$.
\end{enumerate}
\end{lemma}
{\bf Proof:} The symmetry is clear because of the symmetry of $f_i$, $t_i$ and the fact that $s_i$ is the dual of $h_i$. The fact that $\pi_i$ is $D$-parallel comes from the definition of $D$ and equations \eqref{relation1.1} to \eqref{relation2.4}. Indeed, we have for $X,Y\in\Gamma(TM)$,
\beQ
(D_X\pi_i)Y&=&D_X(\pi_iY)-\pi_i(D_XY)\\
&=&D_X(f_iY+h_iY)-\pi_i\left(\nabla_XY+B(X,Y)-\sum_{k=1}^mc_kg(f_kX,Y)\xi_k\right)\\
&=&\nabla_X(f_iY)+B(X,f_iY)-\sum_{k=1}^mc_kg(f_kX,f_iY)\xi_k\\
&&+\nabla_X(h_iY)-A_{h_iY}X-\sum_{k=1}^mc_kg(h_kX,h_iY)\xi_k\\
&&-f_i(\nabla_XY)+h_i(\nabla_XY)-s_iB(X,Y)-t_iB(X,Y)+c_ig(f_iX,Y)\xi_i.
\eeQ
By the use of equations \eqref{relation2.1} and \eqref{relation2.2}, we get
\beQ
(D_X\pi_i)Y&=&c_ig(f_iX,Y)\xi_i-\sum_{k=1}^mc_k\Bigl[g(f_kX,f_iY)+g(h_kX,h_iY)\xi_k\Bigr]\\
&=&c_ig(f_iX,Y)\xi_i-\sum_{k=1}^mc_k\Bigl[g(f_i\circ f_kX+s_i\circ h_kX,Y)\Bigr]\\
&=&0,
\eeQ
since, by \eqref{relation1.1}, we have $f_i\circ f_k+s_i\circ h_k=\delta_i^kf_i$. The computations are analogous for a section $\nu$ of $E$ or for one of the $\xi_i$. \\
The relation $\pi_i\circ\pi_i=\pi_i$ is obvious from the definition of $\pi_i$ and relations \eqref{relation1.1} to \eqref{relation1.4}. Finally, we get immediately that $\sum_{i=1}^m\pi_i=Id_F$ from the definition and assumption \eqref{relation1.0}. 
\hfill$\square$\\ \\
We consider the subsets $F^i$ of $F$ defined by
$$F^i=\left\{\alpha\in F\ |\ \pi_j\alpha=\delta_i^j\alpha\ \text{for any}\ j\in\{1,\cdots,m\}\right\}.$$
Note that, since the $\pi_i$ are symmetric, then the subbundles $F^i$ are orthogonal with respect to $\widetilde{g}$. We finally need a last lemma
\begin{lemma}
For any $i\in\{1,\cdots,m\}$, there exists orthonormal $n_i+1$ parallel sections $\sigma^i_1,\cdots,\sigma^i_{n_i+1}$ of $F^i$.
\end{lemma}
{\bf Proof:} Let $p$ be a point of $M$. For any $i\in\{1,\cdots,m\}$, let $\{v^i_1,\cdots,v^i_{n_1+1}\}$ be  an orthonormal basis of $F_i^p$, which is of dimension $n_i+1$ by the assumption on the rank of $\pi_i$. Moreover, since $\xi_i$ clearly belongs to $F_i$, we can choose $v_i^1=\sqrt{|c_i|}\xi_1(p)$. Thus, we have $\widetilde{g}(v_1^i,v_1^i)={\rm sign}(c_i)$ and $\widetilde{g}(v_k^i,v_k^i)=1$ for $k\in\{2,\cdots,n_i+1\}$. Since the $F^i$ are orthogonal, the set of all $v^k_i$ forms an orthogonal basis of $F_p$. Now, since the connection $D$ is flat and $M$ is simply connected, then for any $i\in\{1,\cdots,m\}$ there exists a family of parallel sections $\sigma^i_1,\cdots,\sigma^i_{n_i+1}$, such that $\sigma^i_k(p)=v^i_k$. Moreover, since $D$ is compatible with the metric $\widetilde{g}$, then the sections are orthonormal. Finally, since the maps $\pi$ are $D$-parallel, then, for any $i\in\{1,\cdots,m\}$ and any $k\in\{1,\cdots,n_i+1\}$, $\pi_i(\sigma^i_k)=\sigma^i_k$, that is $\sigma^i_k$ is a section of $F^i$. This concludes the proof of the lemma.
\hfill$\square$\\ \\
We will construct now the isometric immersion from $M$ into $\widetilde{P}$. For this, we consider the following functions. For $i\in\{1,\cdots,m\}$ and $k\in\{1,\cdots,n_i+1\}$, let $\varphi_k^i$ be defined by
$$\varphi_k^i=\widetilde{g}(\sigma^i_k,\xi_i).$$
The candidate for the isometric immersion is 
$$\varphi:M\lgra \mathbb{E}_1\times\cdots\times\mathbb{E}_m,$$
where $\mathbb{E}_i$ is the Euclidean space $\RR^{n_i+1}$ if $c_i>0$ and the Minkowski space $\mathbb{L}^{n_i+1}$ if $c_i<0$. We will show that the map $\varphi$ goes into $\widetilde{P}\subset \mathbb{E}_1\times\cdots\times\mathbb{E}_m$ and satisfies all the properties stated in Theorem \ref{thm1}.\\
First, we have 
\beQ
\frac{1}{c_i}=\widetilde{g}(\xi_i.\xi_i)&=&\sum_{k=1}^{n_i+1}\widetilde{g}(\xi_i,\sigma^i_k)^2\;\widetilde{g}(\sigma^i_k,\sigma^i_k)\\
&=&{\rm sign}(c_i)\widetilde{g}(\xi_i,\sigma^i_1)^2+\sum_{k=2}^{n_i+1}\widetilde{g}(\xi_i,\sigma^i_k)^2\\
&=&{\rm sign}(c_i)(\varphi_1^i)^2+\sum_{k=2}^{n_i+1}(\varphi_k^i)^2.
\eeQ
Thus, we get that $(\varphi_1^i,\cdots,\varphi_{n_i+1}^i)\in M_i$, the $n_i$-dimensional simply connected space form of curvature $c_i$, and so, $\varphi(M)$ lies in $\widetilde{P}$.\\
Now, we will show that $\varphi$ is an immersion. For this, let $p\in M$ and $v\in T_pM$ so that $\varphi_*(v)=0$. From the definition of $\varphi$, the fact that $\varphi_*(v)=0$ implies
$$\widetilde{g}(\sigma_k^i(p),\pi_iv)=\widetilde{g}(v_k^i,\pi_iv)=0,$$
for any $i\in\{1,\cdots,m\}$ and $k\in\{1,\cdots,n_i+1\}$. Since, for any $i$,  $\{v_k^i\}$ is an orthonormal basis of $F_p^i$, we get that $\pi_iv=0$. Moreover, from Lemma \ref{lempi}, $\sum_{i=1}^m\pi_i=Id_F$, then $v=0$. This holds for any $p$ and any $v$, so we get that $\varphi$ is an immersion.\\
Moreover, for $v,w\in T_pM$, we have
\beQ
\left\langle \varphi_*(v),\varphi_*(w)\right\rangle&=&\sum_{i=1}^m\left({\rm sign}(c_i) \widetilde{g}(\sigma_1^i,\pi_iv)\widetilde{g}(\sigma_1^i,\pi_iw)+\sum_{k=2}^{n_i+1}\widetilde{g}(\sigma_k^i,\pi_iv)\widetilde{g}(\sigma_k^i,\pi_iw)\right)\\
&=&\sum_{i=1}^m\widetilde{g}(\pi_iv,\pi_iw)\\
&=&\widetilde{g}(v,w),
\eeQ
where $\langle\cdot,\cdot\rangle$ is the (pseudo)-Euclidean metric on $\mathbb{E}_1\times\cdots\times\mathbb{E}_m$. Hence, $\varphi$ is an isometric immersion from $M$ into $\widetilde{P}$.\\
Now, we define the following bundle isomorphism 
$$\Phi:F\lgra T(\mathbb{E}_1\times\cdots\times\mathbb{E}_m)_{|\varphi(M)},$$
by $\Phi(\sigma_k^i)=e_k^i$, where $\{e^i_1,\cdots,e^i_{n_i+1}\}$ is the canonical frame of $T\mathbb{E}_i$ restricted to $\varphi(M)$. For $X\in\Gamma(TM)$, we have
\beQ
\Phi(X)&=&\sum_{i=1}^m\sum_{k=1}^{n_i+1}\widetilde{g}(X,\sigma_k^i)e_k^i\\
&=&\sum_{i=1}^m\sum_{k=1}^{n_i+1}\widetilde{g}(\pi_iX,\sigma_k^i)e_k^i\\
&=&\varphi_*(X).
\eeQ
Moreover, for any $i\in\{1,\cdots,m\}$, we have
$$\Phi(\xi_i)=\sum_{k=1}^{n_i+1}\widetilde{g}(\xi_i,\sigma_k^i)e_k^i.=\sum_{k=1}^{n_i+1}\varphi_k^ie_k^i.$$
Hence, $\Phi(\xi_i)$ is the normal direction of $M_i$ in $\mathbb{E}_i$. Since $\Phi$ is an isometry of the fibers, we deduce that $\Phi(E)$ is the normal bundle $T^{\perp}\varphi(M)$ of $\varphi(M)$ in $\widetilde{P}$. We denote by $\widetilde{\varphi}$ the restriction of $\Phi$ to $E$. It is clear that $\widetilde{\varphi}$ is an isomorphism of vector bundles between $E$ and $T^{\perp}\varphi(M)$.\\
Since $\Phi$ sends the orthonormal parallel sections $\{\sigma_k^i\}$ of $F$ onto the orthonormal parallel sections $\{e_k^i\}$ of $\mathbb{E}_1\times\cdots\times\mathbb{E}_m$, we have
\beqt\label{connect1}
\Phi(D_XY)=\nabla^0_{\varphi_*(X)}\varphi_*(Y),
\eeqt
\beqt\label{connect2}
\Phi(D_X\nu)=\nabla^0_{\varphi_*(X)}\widetilde{\varphi}(\nu),
\eeqt
\beqt\label{connect3}
\Phi(D_X\xi_i)=\nabla^0_{\varphi_*(X)}\Phi(\xi_i),
\eeqt
where $\nabla^0$ is the Levi-Civita connection of $\mathbb{E}_1\times\cdots\times\mathbb{E}_m$.\\
For any $i\in\{1,\cdots,m\}$, we define the map $\widetilde{\pi}_i=\Phi\circ\pi_i\circ\Phi^{-1}$. From this definition, it is clear that $\widetilde{\pi}_i(e_k^j)=\delta_i^je_k^j$. Then, it follows that the maps $\widetilde{\pi}_i$ are symmetric, parallel along $\varphi(M)$ and satisfy $\widetilde{\pi}_i\circ\widetilde{\pi}_j=\delta_i^j\widetilde{\pi}_i$ and $\sum_{i=1}^m\widetilde{\pi}_i={\rm Id}_{T(\mathbb{E}_1\times\cdots\times\mathbb{E}_m)}$. Thus, it is clear that these maps are the restrictions on $\varphi(M)$ of the projections on each factor $T\mathbb{E}_i$ of $T(\mathbb{E}_1\times\cdots\times\mathbb{E}_m)$.\\
Moreover, from the definition of $\widetilde{\pi}_i$, we deduce immediately that
 $$\widetilde{\pi}_i(\varphi_*X)=\varphi_*(f_iX)+\Phi(h_iX),$$
and
$$\widetilde{\pi}_i(\Phi\xi)=\varphi_*(s_iX)+\Phi(t_i\xi).$$
Indeed, we have
$$\widetilde{\pi}_i(\varphi_*X)=\Phi(\pi_i(X))=\Phi(f_iX+h_iX)=\varphi_*(f_iX)+\widetilde{\varphi}(h_iX),$$
and
$$\widetilde{\pi}_i(\widetilde{\varphi}(\nu))=\Phi(\pi_i(\nu))=\Phi(s_i\nu+h_i\nu)=\varphi_*(s_i\nu)+\widetilde{\varphi}(t_i\nu).$$
Finally, we will prove that the second fundamental form is given by $B$ and the normal connection is given by $\overline{\nabla}$. From Equation \eqref{connect1}, we have
\beQ
\nabla^0_{\varphi_*(X)}\varphi_*(Y)&=&\Phi(D_XY)\\
&=&\Phi(\nabla_XY+B(X,Y)-\sum_{i=1}^mc_ig(f_iX,Y)\xi_i)\\
&=&\varphi_*(\nabla_XY)+\widetilde{\varphi}(B(X,Y))-\sum_{i=1}^mc_ig(f_iX,Y)\Phi(\xi_i).
\eeQ
Then the normal part in $T\widetilde{P}$ is $\widetilde{\varphi}(B(X,Y))$, which implies that the second fundamental form of the immersion $\varphi$ is $\widetilde{\varphi}\circ B$. Moreover, from Equation \eqref{connect2}, we have
\beQ
\nabla^0_{\varphi_*(X)}\widetilde{\varphi}(\nu)&=&\Phi(D_X\nu)\\
&=&\Phi(\overline{\nabla}_X\nu-A_{\nu}X-\sum_{i=1}^mc_ig(h_iX,\nu)\xi_i)\\
&=&\widetilde{\varphi}(\overline{\nabla}_X\nu)+\varphi_*(A_{\nu}X)-\sum_{i=1}^mc_ig(h_iX,\nu)\Phi(\xi_i).
\eeQ
Thus, the normal part in $T\widetilde{P}$ is $\widetilde{\varphi}(\overline{\nabla}_X\nu)$ and we deduce that $\nabla_X^{\perp}\widetilde{\varphi}(\nu)=\widetilde{\varphi}(\overline{\nabla}_X\nu)$. Then, we get $\nabla^{\perp}\widetilde{\varphi}=\widetilde{\varphi}\overline{\nabla}$. This concludes the proof of the existence in Theorem \ref{thm1}.\\ \\
Now, we will prove the uniqueness of this isometric immersion up to an isometry of $\widetilde{P}$. This follows directly from the following proposition.
\begin{prop}
Let $\varphi,\varphi':M\lgra \widetilde{P}$ be two isometric immersions with respective normal bundles $E,E'$ and second fundamental forms $B,B'$. Let $f_i$, $h_i$ and $f_i'$, $h_i'$ be the $(1,1)$-tensors defined by \eqref{relationfhst} for $\varphi$ and $\varphi'$ respectively. Assume that
\begin{enumerate}[i)]
\item $f_iX=f_i'X$ for any $i\in\{1,\cdots,m\}$ and $X\in\Gamma(TM)$,
\item there exists an isometry of vector bundles $\phi:E\lgra E'$ so that
$$\phi(h_iX)=h_i'X,$$
$$\phi(B(X,Y))=B'(X,Y),$$
$$\phi(\overline{\nabla}_X\nu)=\overline{\nabla}'_X\phi(\nu),$$
\end{enumerate}
for any $i\in\{1,\cdots,m\}$, $X,Y\in\Gamma(TM)$ and $\nu\in E$.\\
Then, there exists an isometry $\alpha$ of $\widetilde{P}$ such that $\varphi'=\alpha\circ\varphi$ and $\alpha_{*|E}=\phi$.
\end{prop}
{\bf Proof:} We give the complete proof for $c_m\neq0$, the case $c_m$ can be proven with a minor modification.\\
As previously, we denote by $\{e_k^i\}$, for $i\in\{1,\cdots,m\}$ and $k\in\{1,\cdots,n_i+1\}$ the canonical frame of $\mathbb{E}_1\times\cdots\times\mathbb{E}_m$. Hence, we denote by $\varphi_k^i$ and $(\varphi')_k^i$ the components of $\varphi$ and $\varphi'$ respectively in the frame $\{e_k^i\}$. We consider the map $G:M\lgra GL(\mathbb{E}_1\times\cdots\times\mathbb{E}_m)$ defined by
$$G_p(\varphi_*(X))=\varphi_*'(X),$$
$$G_p(\nu)=\phi(\nu),$$
$$G_p(\xi_i)=\xi_i',$$
for any $i\in\{1,\cdots,m\}$, $X,Y\in\Gamma(TM)$ and $\nu\in E$, and where $\xi_i$ and $\xi_i'$ are defined by
$$\xi_i=\sum_{k=1}^{n_i+1}\varphi_k^ie_k^i$$
and 
$$\xi_i'=\sum_{k=1}^{n_i+1}(\varphi')_k^ie_k^i.$$
We will show that the map $G$ is constant, that is, that it does not depend on the point $p$. First of all, we remark that for any $i$ and any $X\in\Gamma(TM)$, we have $\nabla^0_X\xi_i=\pi_i(\varphi_*(X)),$ where $\pi_i$ is the projection on $T\mathbb{E}_i$ and $\nabla^0$ is the Levi-Civita connection of $\mathbb{E}_1\times\cdots\times\mathbb{E}_m$. Now, we will show that $\nabla^0G=0$, or equivalently that $\nabla^0_{X}(G(V))-G(\nabla^0_XV)=0$ for any $X\in\Gamma(TM)$ and $V\in\Gamma(T(\mathbb{E}_1\times\cdots\times\mathbb{E}_m)_{|\varphi(M)})$.\\
First, for $V\in\varphi_*(TM)$, that is, $V=\varphi_*(Y)$ with $Y$ tangent to $M$, we have
\beQ
\nabla^0_{X}(G(\varphi_*(Y)))-G(\nabla^0_{X}\varphi_*(Y))&=&\nabla^0_{X}\varphi_*'(Y)-G(\varphi_*(\nabla_XY)+B(X,Y))\\
&&+\sum_{i=1}^mc_ig(f_i(X),Y)\xi_i'\\
&=&\varphi_*'(\nabla_XY)+B'(X,Y)-\sum_{i=1}^mc_ig(f_i'(X),Y)\xi_i'\\
&&-G(\varphi_*(\nabla_XY)+B(X,Y))+\sum_{i=1}^mc_ig(f_i(X),Y)\xi_i'\\
&=&0,
\eeQ
since, by assumption, $f_iX=f_i'X$ and $\phi(B(X,Y))=B'(X,Y)$.
Now, if $\nu\in E$, we have
\beQ
\nabla^0_{X}(G(\nu))-G(\nabla^0_{X}\nu)&=&-A'_{\phi(\nu)}X+\overline{\nabla}'_X\phi(\nu)-\sum_{i=1}^mc_i\widetilde{g}(\phi(\nu),h_i'(X))\xi'\\
&&-G\left( -A_{\nu}X+\overline{\nabla}_X\nu-\sum_{i=1}^mc_i\widetilde{g}(\phi(\nu),h_i(X))\xi\right)\\
&=&0, 
\eeQ
since by assumption, $\phi(B(X,Y))=B'(X,Y)$, $\phi(h_iX)=h_i'X$ and $\phi(\overline{\nabla}_X\nu)=\overline{\nabla}'_X\phi(\nu)$.\\
Finally, we have
\beQ
\nabla^0_{X}(G(\xi))-G(\nabla^0_{X}\xi)&=&\pi_i(\varphi_*'(X))-G(\pi_i(\varphi_*(X)))\\
&=&\varphi_*'(f_i'X)+h_i'X-G(\varphi_*(f_iX)+h_iX)\\
&=&0,
\eeQ
since $f_iX=f_i'X$ and $\phi(h_iX)=h_i'X$.\\
Hence, we get that the map $G$ is constant along $M$.
\hfill$\square$\\ \\
Uniqueness up to rigid motion is proved, which concludes the proof of Theorem \ref{thm1}.
\hfill$\square$
\section{Associated families of minimal surfaces and pluriminimal K\"alher hypersufaces}
In this section, we use Theorem \ref{thm1} to prove the existence of  associated families of minimal surfaces into the multiproduct $\widetilde{P}$.\\
Let $(\Sigma,g)$ be an oriented Riemannian surface. We denote by $J$ its complex structure, that is, the rotation of angle $\frac{\pi}{2}$ on $TM$. For any $\theta\in\RR$, we set $\mathcal{R}_{\theta}=\cos(\theta)I+\sin(\theta)J$. Remark, that $\mathcal{R}_{\theta}$ is parallel. First, we have the following proposition.
\begin{prop}\label{prop1}
Assume that $(\Sigma,g,E,\overline{g},\overline{\nabla},B,f_i,h_i,t_i)$ satisfies the compatibility equation for $\widetilde{P}$ and that $B$ is trace-free for any $\xi\in E$, then $(\Sigma,g,E,\overline{g},\overline{\nabla},B_{\theta},f_{i,\theta},h_{i,\theta},t_{,\theta})$ also satisfies the compatibility equations for $\widetilde{P}$, where
\begin{align}
&B_{\theta}(X,Y)=B(\mathcal{R}_{\theta}X,Y),&\nonumber\\
&f_{i,\theta}=\mathcal{R}_{\theta}\circ f_i\circ\mathcal{R}^{-1}_{\theta},&\nonumber\\
&h_{i,\theta}=h_i\circ\mathcal{R}^{-1}_{\theta},&\nonumber\\
&t_{i,\theta}=t_i.&\nonumber
\end{align}
Moreover, $B_{\theta}$ is also trace-free for any $\xi\in E$.
\end{prop}
{\bf Proof:} First, from the definition of $f_{i,\theta}$, $h_{i,\theta}$ and $t_{i,\theta}$ and the fact that 
\begin{align}
&f_i\circ f_j+s_i\circ h_j=\delta_i^jf_i,&\nonumber\\
&t_i\circ t_j+h_i\circ s_j=\delta_i^jt_i,&\nonumber \\
&f_i\circ s_j+s_i\circ t_j=\delta_i^js_i,&\nonumber\\
&h_i\circ f_j+t_i\circ h_j=\delta_i^jh_i,&\nonumber
\end{align} we get immediately that $f_{i,\theta}$ and $t_{i,\theta}$ are symmetric and
\begin{align}
&f_{i,\theta}\circ f_{j,\theta}+s_{i,\theta}\circ h_{j,\theta}=\delta_i^jf_{i,\theta},&\nonumber\\
&t_{i,\theta}\circ t_{j,\theta}+h_{i,\theta}\circ s_{j,\theta}=\delta_i^jt_{i,\theta},&\nonumber \\
&f_{i,\theta}\circ s_{j,\theta}+s_{i,\theta}\circ t_{j,\theta}=\delta_i^js_{i,\theta},&\nonumber\\
&h_{i,\theta}\circ f_{j,\theta}+t_{i,\theta}\circ h_{j,\theta}=\delta_i^jh_{i,\theta}.&\nonumber
\end{align}
It is also clear that with this definition, the rank of $\pi_{i,\theta}$ is the same that the rank of $\pi_i$ and that 
$$\sum_{i=1}^mf_{i,\theta}=Id_{TM},\quad \sum_{i=1}^mt_{i,\theta}=Id_{E}\quad\text{and}\quad \sum_{i=1}^mh_{i,\theta}=0.$$

 Now, we will show that analogues of Equations \eqref{relation1.1}-\eqref{relation1.3} are satisfied for $(\Sigma,g,E,\overline{g},\overline{\nabla},B_{\theta},f_{i,\theta},h_{i,\theta},t_{,\theta})$. First, we have for $X,Y$ tangent to $\Sigma$
\beQ
\nabla_X f_{i,\theta}Y-f_{i,\theta}(\nabla_XY)&=&\nabla_X(\mathcal{R}_{\theta}f_i\mathcal{R}_{\theta}^{-1}Y)-\mathcal{R}_{\theta}f_i\mathcal{R}_{\theta}^{-1}(\nabla_XY)\\
&=&\mathcal{R}_{\theta}\nabla_X(f_i\mathcal{R}_{\theta}^{-1}Y)-\mathcal{R}_{\theta}f_i\mathcal{R}_{\theta}^{-1}(\nabla_XY),
\eeQ
since $\mathcal{R}_{\theta}$ is parallel. Moreover, using \eqref{relation1.1}, we get
\beQ
\nabla_X f_{i,\theta}Y-f_{i,\theta}(\nabla_XY)&=&\mathcal{R}_{\theta}A_{h_i(\mathcal{R}_{\theta}^{-1}Y)}X+\mathcal{R}_{\theta}s_i\left(B(X,\mathcal{R}_{\theta}^{-1}Y)\right)\\&&
+ \mathcal{R}_{\theta}f_i\mathcal{R}_{\theta}^{-1}(\nabla_XY)-\mathcal{R}_{\theta}f_i(\nabla_X\mathcal{R}_{\theta}^{-1}Y)\\
&=&A^{\theta}_{h_{i,\theta}Y}X+s_{i,\theta}\left(B_{\theta}(X,Y)\right),
\eeQ
which is the desired equation. The two other equations can be shown in a similar way.\\
Finally, we prove that Gauss, Codazzi and Ricci equations are also fulfilled.\\
First we consider the Gauss equation. We notice that, for a surface, we have
\begin{eqnarray*}\label{gauss}
&&\sum_{i=1}^mc_i\bigg[ \left\langle f_{i,\theta}Y,Z\right\rangle f_{i,\theta}X-\left\langle f_{i,\theta}X,Z\right\rangle f_{i,\theta}Y \bigg]+{A_{\theta}}_{B_{\theta}(Y,Z)}X-{A_{\theta}}_{B_{\theta}(X,Z)}Y\\
\\
%&&=\sum_{i=1}^mc_i\bigg[ \left\langle %\mathcal{R}_{\theta}f_i\mathcal{R}_{\theta}^{-1}Y,Z\right\rangle\mathcal{R}_{\theta}f_i\mathcal{R}_{\theta}^{-1}X-\left\langle %\mathcal{R}_{\theta}f_i\mathcal{R}_{\theta}^{-1}X,Z\right\rangle\mathcal{R}_{\theta}f_i\mathcal{R}_{\theta}^{-1}Y %\bigg]+{A_{\theta}}_{B_{\theta}(Y,Z)}X-{A_{\theta}}_{B_{\theta}(X,Z)}Y\\
&&=\sum_{i=1}^mc_i\bigg[  \mathcal{R}_{\theta}f_i\mathcal{R}_{\theta}^{-1}X\wedge\mathcal{R}_{\theta}f_i\mathcal{R}_{\theta}^{-1}Y \bigg]Z+\det A_{\theta}\\
&&=\sum_{i=1}^mc_i\det \mathcal{R}_{\theta}f_i\mathcal{R}_{\theta}^{-1}+\det A_{\theta}=\sum_{i=1}^mc_i\det f_i+\det A=R(X,Y)Z
\end{eqnarray*}
since determinants are invariant under rotations. Hence Gauss equation is satisfied.\\
Let $\widetilde{\nabla}_XA^{\nu}_{\theta}=\nabla_XA^{\nu}_{\theta}Y-A^{\nu}_{\theta}\nabla_XY-{A^{\nu}_{\theta}}_{\nabla^{\perp}_{X}\nu}Y$. Considering now Codazzi equation, we have, using the property of $h_i$, 
\begin{eqnarray*}
(\widetilde{\nabla}_XA^{\nu}_{\theta})Y-(\widetilde{\nabla}_YA^{\nu}_{\theta})X
&=&\mathcal{R}_{\theta}\Big[(\widetilde{\nabla}_XA^{\nu})Y-(\widetilde{\nabla}_YA^{\nu})X\Big]\\
&=&\mathcal{R}_{\theta}\sum_{i=1}^mc_i\bigg[ f_{i}Y\left\langle X,s_{i}\nu\right\rangle-f_{i}X\left\langle  Y,s_{i}\nu\right\rangle \bigg]\\
&=&\mathcal{R}_{\theta}\sum_{i=1}^mc_i f_i(X\wedge Y)s_i\nu=\sum_{i=1}^mc_i \mathcal{R}_{\theta}f_i(X\wedge Y)\mathcal{R}^{-1}_{\theta}\mathcal{R}_{\theta}s_i\nu
\\&=&\sum_{i=1}^mc_i (\mathcal{R}_{\theta}f_i\mathcal{R}^{-1}_{\theta}X\wedge Y)s_{i,\theta}\nu=
\sum_{i=1}^mc_i (f_{i,\theta}X\wedge Y)s_{i,\theta}\nu
\\
&=&\sum_{i=1}^mc_i\bigg[ f_{i,\theta}Y\left\langle X,s_{i,\theta}\nu\right\rangle-f_{i,\theta}X\left\langle  Y,s_{i,\theta}\nu\right\rangle \bigg]=\sum_{i=1}^mc_i\bigg[ f_{i,\theta}Y\left\langle h_{i,\theta} X,\nu\right\rangle-f_{i,\theta}X\left\langle h_{i,\theta} Y,\nu\right\rangle \bigg],
\end{eqnarray*}
and Codazzi is satisfied.\\
Similarly we get for the Ricci equation, using the properties of the wedge product
\begin{eqnarray*}\label{ricci}
R^{\perp}(X,Y)\xi&=&\sum_{i=1}^mc_i\bigg[  h_i\mathcal{R}_{\theta}^{-1}X\wedge h_i\mathcal{R}_{\theta}^{-1}Y \bigg]\xi+B_{\theta}({A_{\theta}}_{\xi}Y,X)-B_{\theta}({A_{\theta}}_{\xi}X,Y)\\
&=&\sum_{i=1}^mc_i\bigg[  h_iX\wedge h_iY \bigg]\xi+B({A}_{\xi}Y,X)-B({A}_{\xi}X,Y)
\end{eqnarray*}
Since the surface is minimal, the shape operator anti-commutes with $J$ and we have indeed  
\begin{eqnarray*}
B^{\nu}_{\theta}({A_{\theta}}_{\xi}Y,X)-B^{\nu}_{\theta}({A_{\theta}}_{\xi}X,Y)&=&\langle[{A_{\theta}}_{\nu},{A_{\theta}}_{\xi}]X,Y\rangle=\langle (\mathcal{R}_{\theta}A_{\nu} \mathcal{R}_{\theta}A_{\xi}-\mathcal{R}_{\theta}A_{\xi}\mathcal{R}_{\theta}A_{\nu})X,Y\rangle\\
&=&
\langle (A_{\nu} \mathcal{R}_{\theta}^{-1}\mathcal{R}_{\theta}A_{\xi}-A_{\xi}\mathcal{R}_{\theta}^{-1}\mathcal{R}_{\theta}A_{\nu})X,Y\rangle=\langle[A_{\nu},A_{\xi}]X,Y\rangle.
\end{eqnarray*}
Finally, let $(e_1,e_2=Je_1)$ be a local orthonromal frame of $\Sigma$. We have 
\begin{eqnarray*}
\tr(B_{\theta})=B_{\theta}(e_1,e_1)+B_{\theta}(e_2,e_2)&=&B(\mathcal{R}_{\theta}e_1,e_1)+B(\mathcal{R}_{\theta}e_2,e_2)\\
&=& \cos\theta\left[B(e_1,e_1)+B(e_2,e_2)\right]=0,
\end{eqnarray*}
since $B$ is trace-free.
\hfill$\square$\\ \\ 
From this proposition, we can prove easily the following theorem about associated families of minimal surfaces in multiproducts. Namely, we get the following statement.
\begin{thm}\label{thm2}
Let $\Sigma$ be a simply connected surface and $x:M\lgra \widetilde{P}$ be a conformal minimal immersion with normal bundle $E$, second fundamental form $B$ and normal connection $\nabla^{\perp}$. Let $f_i$, $h_i$, $s_i$ and $t_i$ be the $(1,1)$-tensors induced by the projections $\pi_i$. Let $p_0\in\Sigma$. Then, there exists a unique family $(x_{\theta})_{\theta\in\RR}$ of conformal minimal immersions $x_{\theta}:\Sigma\lgra\widetilde{P}$ so that
\begin{enumerate}[i)]
\item $x_{\theta}(p_0)=x(p_0)$ and $d({x_{\theta}})_{p_0}=(dx)_{p_0}$,
\item the metric induced by $X$ and $X_{\theta}$ are the same,
\item the second fundamental form fo $x_{\theta}(\Sigma)$ in $\widetilde{P}$ is given by $B_{\theta}(X,Y)=B(R_{\theta}X,R_{\theta}Y)$, for any $X,Y\in\Gamma(T\Sigma)$.
\item for any $i\in\{1,\cdots,m\}$, $X\in\Gamma(T\Sigma)$ and $\xi\in\Gamma(E)$,
$$\pi_i(dx_{\theta}X)=dx_{\theta}(f_{i,\theta}X)+h_{i,\theta}X\quad \text{and}\quad \pi_i(\xi)=dx_{\theta}(s_{i,\theta}X)+t_{i,\theta}X,$$
\end{enumerate}
Moreover, $x_0=x$ and the family $(x_{\theta})_{\theta\in\RR}$ is continuous with respect to $\theta$.
\end{thm}

{\bf Proof:} We just proved that $(\Sigma,g,E,\overline{g},\overline{\nabla},B_{\theta},f_{i,\theta},h_{i,\theta},t_{,\theta})$ satisfies the compatibility equations for each $\theta$. The theorem is then a direct consequence of theorem \ref{thm1}. The continuity is ensured by the construction of Theorem \ref{thm1}.
\hfill$\square$

\section{Surfaces in $\mathbb{S}^2\times\mathbb{S}^2$}
 Let $J$ be the complex structure on $\mathbb{S}^2$. We consider the following complex structures on $\mathbb{S}^2\times\mathbb{S}^2$
\[J_1=(J,\,J),\quad J_2=(J,\,-J).\]
Obviously $J_1$ and $J_2$ commute with each other and the projection $\pi_1$ and $\pi_2$ on each of the factors are given by 
\[\pi_1=\frac{\id+J_1J_2}{2},\quad \pi_2=\frac{\id-J_1J_2}{2}.\]
From equation \eqref{extcurvature} we get  
\begin{eqnarray}\label{curvJ}
\widetilde{R}(X,Y)Z=\bigl[\langle Y,Z\rangle X-\langle X,Z\rangle Y\bigr]+\bigl[\langle J_1 Y,J_2Z\rangle J_1J_2X-\langle J_1X,J_2Z\rangle J_1J_2Y\bigr].\end{eqnarray}
Let now $\Sigma$ be a surface isometrically immersed into $\mathbb{S}^2\times\mathbb{S}^2$. For $i\in\{1,2\}$, we define four operators $j_i:T\Sigma\rightarrow T\Sigma$,  $k_i:T\Sigma\rightarrow N\Sigma$ ,$ l_i:N\Sigma\rightarrow T\Sigma$ and $m_i:N\Sigma\rightarrow N\Sigma$ such that $J_i=j_i+k_i+l_i+m_i$.\\
From $J_1J_2=J_2J_1$ we get the following equations
\begin{eqnarray}
&&j_1j_2+l _1k_2=j_2j_1+l_2k_1,\label{relation5com1}\\
&&k_1j_2+m_1k_2=k_2j_1+m_2k_1,\label{relation5com2}\\
&&j_1l_2+l_1m_2=j_2l_1+l_2m_1,\label{relation5com3}\\
&&k_1l_2+m_1m_2=k_2l_1+m_2m_1.\label{relation5com4}
\end{eqnarray}
The property $J_i^2=-\Id$ yields
\begin{eqnarray}
&&j_i^2+l _ik_i=-\Id_{T\Sigma},\label{relation5cx1}\\
&&k_ij_i+m_ik_i=0,\label{relation5cx2}\\
&&j_il_i+l_im_i=0,\label{relation5cx3}\\
&&k_il_i+m_i^2=-\Id_{N\Sigma}.\label{relation5cx4}
\end{eqnarray}
Moreover the fact that the operators $J_i$ are antisymmetric implies the antisymmetry of the operators $j_i$ as well as the property $\langle k_iX,\nu\rangle=-\langle X,l_i\nu\rangle$.\\
The parallelity of $J_i$ gives
\begin{align}
&\nabla_X(j_iY)-j_i(\nabla_XY)=A_{k_iY}X+l_i(B(X,Y)),& \label{relation5.1}\\
&\nabla^{\perp}_X(k_iY)-k_i(\nabla_XY)=m_i(B(X,Y))-B(X,j_iY),& \label{relation5.2}\\
&\nabla_X^{\perp}(m_i\xi)-m_i(\nabla^{\perp}_X\xi)=-B(l_i\xi,X)-k_i(A_{\xi}X),& \label{relation5.3}\\
&\nabla_X(l_i\xi)-l_i(\nabla^{\perp}_X\xi)=-j_i(A_{\xi}X)+A_{m_i\xi}X.&\label{relation5.4}
\end{align}
Finally, from \eqref{curvJ}, we get the Gauss equation
\begin{eqnarray}\label{gaussex}
K&=&\frac12\bigg[ 1+ \Big(\langle j_1e_1,j_2e_2\rangle+\langle k_1e_1,k_2e_2\rangle\Big)\Big(\langle j_1e_2,j_2e_1\rangle+\langle k_1e_2,k_2e_1\rangle\Big)\\ \nonumber
&& - \Big(\langle j_1e_1,j_2e_1\rangle+\langle k_1e_1,k_2e_1\rangle\Big)\Big(\langle j_1e_2,j_2e_2\rangle-\langle k_1e_2,k_2e_2\rangle\Big)\bigg]\\ \nonumber
&&+2|H|^2-\frac{|B|^2}{2},
\end{eqnarray}
the Codazzi equation
\begin{eqnarray}\label{codazziex}\nonumber
&&(\nabla_XB)(Y,Z)-(\nabla_YB)(X,Z)=\\
&&\hspace{1cm}\frac12\bigg[\langle Y,(j_1j_2+l_1k_2)Z\rangle (k_1j_2+m_1k_2)X-\langle X,(j_1j_2+l_1k_2)Z\rangle (k_1j_2+m_1k_2)Y\bigg], 
\end{eqnarray}
and the Ricci equation
\begin{eqnarray}\label{ricciex}
K^{\perp}&=&\bigg[ \left\langle k_1j_2+m_1k_2)e_2,\nu_1\right\rangle \left\langle(k_1j_2+m_1k_2)e_1,\nu_2\right\rangle \\ \nonumber
&&\left\langle k_1j_2+m_1k_2)e_1,\nu_1\right\rangle \left\langle(k_1j_2+m_1k_2)e_2,\nu_2\right\rangle\bigg]\\ \nonumber
&&+\langle [A_{\nu_2},A_{\nu_1}]e_1,e_2]\rangle.
\end{eqnarray} 
\begin{remark}
In \cite{TU} the Gauss, Codazzi and Ricci equations are expressed with the help of the two K\"ahler functions $C_1$ and $C_2:\Sigma\rightarrow \mathbb{R}$ defined by $\varphi^*\omega_i=C_i\omega_{\Sigma},\,i=1,2$, with $\omega_{\Sigma}$ the area form on $\Sigma$. A tidious but straightforward computation shows that those two formulations are equivalent.
\end{remark}
Now, we are able to reformulate the main theorem in the case of $\mathbb{S}^2\times\mathbb{S}^2$ in terms of complex structures instead of projections on each factor.
\begin{cor}
Let $(\Sigma^2,g)$ be a Riemannian surface and $(E,\langle\cdot,\cdot\rangle_E,\nabla^E)$ a rank $2$ vector bundle over $\Sigma$ endowed with a scalar product and a compatible connection. Suppose that there exists a symmetric $(2,1)$-tensor field $B: T\Sigma \times T\Sigma\rightarrow E$ and eight operators $j_i:T\Sigma\rightarrow T\Sigma$,  $k_i:T\Sigma\rightarrow E$, $ l_i:E\rightarrow T\Sigma$ and $m_i:E\rightarrow E$, $i=1,2$ satisfying conditions \eqref{relation5com1} to \eqref{relation5.4} and the Gauss, Codazzi and Ricci equations \eqref{gaussex}, \eqref{codazziex} and \eqref{ricciex}. Then, there exists a unique (up to isometries of \,$\mathbb{S}^2\times\mathbb{S}^2$) isometric immersion from $\Sigma$ into $\mathbb{S}^2\times\mathbb{S}^2$ with $E$ as normal bundle, $B$ as second fundamental form and such that the restrictions of the complex structures $J_i$ over $\Sigma$ are given by $j_i$, $k_i$, $l_i$ and $m_i$.
\end{cor}
{\bf Proof:} Define the following operators
\begin{eqnarray*}
f_1&=&\Id_{T\Sigma}+j_1j_2+l_1k_2,\quad f_2=\Id_{T\Sigma}-j_1j_2-l_1k_2\\
h_1&=&k_1j_2+m_1k_2,=-h_2,\\
s_1&=&j_1l_2+l_1m_2=-s_2 \\\
t_1&=&\Id_E+k_1l_2+m_1m_2,\quad t_2=\Id_E-k_1l_2-m_1m_2.
\end{eqnarray*} 
We can show easily that equations \eqref{relation5com1} to \eqref{ricciex} imply that these operators satisfy the compatibility equations for  $\mathbb{S}^2\times\mathbb{S}^2$ given by Definition \ref{defcomp}. The conclusion follows easily from Theorem \ref{thm1}. \qed
\begin{remark}
We remind (see for example \cite{TU}) that an immersion $\varphi:\Sigma\longrightarrow\mathbb{S}^2\times\mathbb{S}^2$ is called complex if it is complex with respect to $J_1$ or $J_2$. It is called Lagrangian, if it is Lagrangian with respect to $J_1$ or $J_2$.
\paragraph{First Case:} $\Sigma$ is a complex surface with respect to one of the complex structures $J_i$. Then it is automatically minimal. Moreover $k_i=l_i=0$, $j_i$ and $m_i$ are parallel complex structures on $T\Sigma$ and $N\Sigma$  respectively, and $j_1$ commutes with $j_2$, as well as $m_1$ with $m_2$. Assume without loss of generality that $\Sigma$ is complex with respect to $J_1$, then the Gauss, Codazzi and Ricci equations simplify to 
\begin{eqnarray*}
&&K=\frac12\bigg[ 1+ \langle j_1e_1,j_2e_2\rangle\langle j_1e_2,j_2e_1\rangle- \langle j_1e_1,j_2e_1\rangle\langle j_1e_2,j_2e_2\rangle\bigg]-\frac{|B|^2}{2},\\
&&(\nabla_XB)(Y,Z)-(\nabla_YB)(X,Z)=\frac12\bigg[\langle Y,j_1j_2Z\rangle m_1k_2 X-\langle X,j_1j_2Z\rangle m_1k_2Y\bigg], \\
&&K^{\perp}=\bigg[ \left\langle m_1k_2e_2,\nu_1\right\rangle \left\langle m_1k_2e_1,\nu_2\right\rangle-\left\langle m_1k_2 e_1,\nu_1\right\rangle \left\langle m_1k_2e_2,\nu_2\right\rangle\bigg]+\langle [A_{\nu_2},A_{\nu_1}]e_1,e_2]\rangle.
\end{eqnarray*}
Notice that the only examples of complex surfaces with respect to both complex structures $J_1$ and $J_2$ are slices $\mathbb{S}^2\times{p}=\{(x,p)\in\mathbb{S}^2\times\mathbb{S}^2|x\in\mathbb{S}^2\}$ and 
$p\times\mathbb{S}^2=\{(p,x)\in\mathbb{S}^2\times\mathbb{S}^2|x\in\mathbb{S}^2\}$.
\paragraph{Second Case:} $\Sigma$ is Lagrangian with respect to $J_i$, then $j_i=m_i=0$. Assuming again without loss of generality that $\Sigma$ is Lagrangian with respect to $J_1$, the Gauss, Codazzi and Ricci equations simplify in the following way
\begin{eqnarray*}
&&K=\frac12\bigg[ 1- \langle k_1e_1,k_2e_2\rangle\langle k_1e_2,k_2e_1\rangle+ \langle k_1e_1,k_2e_1\rangle\langle k_1e_2,k_2e_2\rangle\bigg] +2|H|^2-\frac{|B|^2}{2},\\
&&(\nabla_XB)(Y,Z)-(\nabla_YB)(X,Z)=\frac12\bigg[\langle Y,l_1k_2)Z\rangle k_1j_2X-\langle X,l_1k_2Z\rangle k_1j_2Y\bigg],\\
&&K^{\perp}=\bigg[ \left\langle k_1j_2e_2,\nu_1\right\rangle \left\langle k_1j_2e_1,\nu_2\right\rangle-\left\langle k_1j_2e_1,\nu_1\right\rangle \left\langle k_1j_2e_2,\nu_2\right\rangle\bigg]+\langle [A_{\nu_2},A_{\nu_1}]e_1,e_2]\rangle.
\end{eqnarray*}
Notice that $\Sigma$ is Lagrangian for both $J_1$ and $J_2$ if and only if it is the product $\varphi(s,t)=(\alpha(s), \beta(t))$ of two curves $\alpha, \beta$ in $\mathbb{S}^2$. The Clifford torus is the only example of a minimal such surface.
\paragraph{Third Case:} $\Sigma$ is Lagrangian with respect to $J_1$ (hence $j_1=m_1=0$) and complex with respect to $J_2$ (hence $k_2=l_2=0$).
\begin{eqnarray*}
&&K=\frac12-\frac{|B|^2}{2},\\
&&(\nabla_XB)(Y,Z)-(\nabla_YB)(X,Z)=0,\\
&&K^{\perp}=\bigg[ \left\langle k_1j_2e_2,\nu_1\right\rangle \left\langle k_1j_2e_1,\nu_2\right\rangle-\left\langle k_1j_2e_1,\nu_1\right\rangle \left\langle k_1j_2e_2,\nu_2\right\rangle\bigg]+\langle [A_{\nu_2},A_{\nu_1}]e_1,e_2]\rangle.
\end{eqnarray*}
The only example of such a surface is the diagonal $\mathbf{D}=\{(x,x)\in\mathbb{S}^2\times\mathbb{S}^2|x\in\mathbb{S}^2\}$.
\end{remark}

\end{document}